\pdfoutput=1
\RequirePackage{ifpdf}
\ifpdf 
\documentclass[pdftex]{sigma}
\else
\documentclass{sigma}
\fi

\usepackage{mathtools}

\usepackage{mathrsfs}

\newcommand{\R}{\ensuremath{\mathbb R}}

\DeclareMathOperator{\tr}{tr}

\let\div\undefined
\newcommand{\div}{\ensuremath{\mathrm{div}}} 
\newcommand{\lapl}{\ensuremath{\Delta}} 

\DeclareMathOperator{\II}{I\!I} 
\newcommand{\Scal}{\ensuremath{R}} 
\DeclareMathOperator{\Ric}{Ric} 
\DeclareMathOperator{\vol}{vol} 

\newcommand{\M}{\ensuremath{\mathscr M}} 

\newcommand{\met}{\ensuremath{\mathscr R}} 

\newcommand{\e}{\varepsilon}
\newcommand{\be}{\begin{equation}}
\newcommand{\ee}{\end{equation}}

\newtheorem{Theorem}{Theorem}[section]
\newtheorem{Corollary}[Theorem]{Corollary}
\newtheorem{Lemma}[Theorem]{Lemma}
\newtheorem{Proposition}[Theorem]{Proposition}
{ \theoremstyle{definition}
\newtheorem{Definition}[Theorem]{Definition}
\newtheorem{question}[Theorem]{Question}

\newtheorem{Remark}[Theorem]{Remark} }

\numberwithin{equation}{section}

\begin{document}

\allowdisplaybreaks

\renewcommand{\thefootnote}{}

\newcommand{\arXivNumber}{2306.17760}

\renewcommand{\PaperNumber}{014}

\FirstPageHeading

\ShortArticleName{A Note about Isotopy and Concordance of Positive Scalar Curvature Metrics}

\ArticleName{A Note about Isotopy and Concordance\\ of Positive Scalar Curvature Metrics\\ on Compact Manifolds with Boundary\footnote{This paper is a~contribution to the Special Issue on Differential Geometry Inspired by Mathematical Physics in honor of Jean-Pierre Bourguignon for his 75th birthday. The~full collection is available at \href{https://www.emis.de/journals/SIGMA/Bourguignon.html}{https://www.emis.de/journals/SIGMA/Bourguignon.html}}}

\Author{Alessandro CARLOTTO~$^{\rm a}$ and Chao LI~$^{\rm b}$}

\AuthorNameForHeading{A.~Carlotto and C.~Li}

\Address{$^{\rm a)}$~Universit\`a di Trento, Dipartimento di Matematica, via Sommarive 14, 38123 Trento, Italy}
\EmailD{\href{mailto:alessandro.carlotto@unitn.it}{alessandro.carlotto@unitn.it}}

\Address{$^{\rm b)}$~New York University - Courant Institute of Mathematical Sciences,\\
\hphantom{$^{\rm b)}$}~251 Mercer Street, New York, NY 10012, USA}
\EmailD{\href{mailto:chaoli@nyu.edu}{chaoli@nyu.edu}}

\ArticleDates{Received July 03, 2023, in final form January 31, 2024; Published online February 13, 2024}

\Abstract{We study notions of isotopy and concordance for Riemannian metrics on manifolds with boundary and, in particular, we introduce two variants of the concept of minimal concordance, the weaker one naturally arising when considering certain spaces of metrics defined by a suitable spectral ``stability'' condition. We develop some basic tools and obtain a rather complete picture in the case of surfaces.}


\Classification{53C21; 53A10}

\begin{flushright}
\begin{minipage}{65mm}
\it Dedicated to Jean-Pierre Bourguignon,\\ with gratitude and admiration, \\ on the occasion of his 75th birthday.
\end{minipage}
\end{flushright}

\renewcommand{\thefootnote}{\arabic{footnote}}
\setcounter{footnote}{0}

 Jean-Pierre Bourguignon has been, for more than half a century, a prominent figure on the mathematical scene. A student of Marcel Berger, his work combines the austerity and strive for generality that characterized the Bourbaki movement (as we can see, for instance, in the monumental treatise \cite{Bes87} on Einstein manifolds or in the earlier monograph \cite{Bes78}, but even in the style of some of his most celebrated research articles such as \cite{Bou81}) with an authentic interest for physical phenomena and their rigorous description. Just to name a few towering contributions in that direction, he has studied, over the years, central themes in fluid dynamics (see, e.g., \cite{BouBre74} with Brezis), Yang--Mills fields (as in~\cite{BouLaw81} with Lawson), spinors and Dirac operators (as in \cite{BouGau92} with Gauduchon) \dots\

 On top of that, Bourguignon has been a strenuous, tireless advocate for mathematics, serving in a variety of managing and institutional roles. Among them, he was for two decades the director of the \emph{Institut des Hautes \'Etudes Scientifiques} in Paris, keeping its high, uncompromised scientific standards in the post-Grothendieck era, before becoming, in 2014, chair of the \emph{European Research Council}. In that position he was successful in defending the need of adequately supporting scientific research with large-scale grants, aimed at funding ambitious projects in all scientific fields, beyond the self-evident needs of experimental disciplines. There is no doubt that such initiatives helped retain in Europe a number of early-career scientists that would have otherwise opted for different paths, thereby forming, in turn, new generations of researchers. For all these reasons, and for his friendly presence over the years, it is a pleasure for us to dedicate him this manuscript, with sincere admiration.
						
\section{Introduction} \label{sec:intro}

One of the recurring themes in Bourguignon's mathematical work is the study of \emph{scalar curvature}. In recent years, this field has flourished, with both significant advances on some classical questions and the opening of new avenues (we refer interested readers to the survey \cite{Car21} and references therein). Among emerging trends, one can certainly include the study of positive scalar curvature metrics on compact manifolds with boundary: in the last decade the emphasis has been shifted from collar-type (i.e., product) boundary conditions to local ones, such as (typically) those defined by a binary constraint involving the mean curvature of the boundary in question, or by the more restrictive requirements one can impose on the second fundamental form.

In short, and oversimplifying things to the extreme, we have learnt (see \cite[Section 3]{CarLi19} and~\cite[Section~4]{BarHan20}) -- among other things -- that the study of spaces of positive scalar curvature (henceforth: PSC) metrics with mean-convex or minimal boundary conditions can be reduced (in the sense of weak homotopy equivalence) to the study of spaces of PSC metrics with \emph{doubling} boundary condition over the same manifold, hence to spaces of PSC metrics with a reflectional symmetry on closed manifolds, which inspired the definition of \emph{reflexive manifolds} (see \mbox{\cite{CarLi19, CarLi21}}). When the compact background manifold has dimension two or three, then these spaces of metrics, when not empty, are actually \emph{contractible} (see, respectively, \cite{CarWu21} and \cite{CarLi21}, the latter crucially building on deep work by Bamler--Kleiner \cite{BamKle19}). On the other hand, we note that such spaces are \emph{not}, in general, homotopically equivalent to the space of PSC metrics with collar boundary over the same background manifold (see, e.g., Remark 41 in \cite{BarHan20} for a simple yet enlightening example, and the main results in \cite{RosWei22} providing a thorough investigation in the high-dimensional regime). The reader may want to compare these results (and, in fact, the discussion we are about to present in this paper) with earlier works by Akutagawa--Botvinnik~\cite{AkutagawaBotvinnik2002conformal,AkutagawaBotvinnik2002relative} and Akutagawa \cite{Akutagawa2003} for different yet partly related contributions, connected to the notion of relative Yamabe invariant.

 That being said, and getting back -- at least momentarily -- to closed manifolds, another interesting class of spaces of Riemannian metrics has clearly emerged in the influential work~\cite{ManSch15} by Mantoulidis and Schoen, which in a way stems from classical ideas and methods dating back to Schoen and Yau (see, in particular, \cite{SchYau79}). Motivated by the basic problem of computing the Bartnik quasi-local mass \cite{Bar89} for ample classes of examples, they studied the space of metrics corresponding (extrinsically speaking) to those induced on (closed) \emph{stable, minimal} hypersurfaces inside an ambient manifold of positive scalar curvature. They ultimately exploited the flexibility properties such a space enjoys to
 prove a fairly definitive result in the so-called minimal case (namely for Bartnik triples with $H\equiv 0$): for any such triple, the value of the Bartnik quasi-local mass equals the lower bound prescribed by the Riemannian Penrose inequality (see \cite{Bra01,HuiIlm01} as well as the recent, striking proof in \cite{AMMO22}) and is actually \emph{never} achieved, except when the data correspond to a round sphere. Their construction, which also crucially contains a smoothing procedure (albeit in spherical symmetry), has then been extensively studied and generalized in different directions; the reader is refereed to the beautiful survey \cite{CabCed21} for a broad-spectrum account on this matter.

Our purpose here is to study analogous phenomena in the aforementioned case of \emph{manifolds with boundary}. To some extent we shall analyze how to properly extend the results in \cite{ManSch15} to that setting. We shall be concerned with some natural questions regarding the equivalence relations of \emph{isotopy} and \emph{concordance}; in particular, we will see how a weak notion of \emph{minimal concordance} naturally arises when considering certain spaces of metrics defined by the corresponding stability condition, but now in the category of free boundary minimal surfaces. The key definitions and some preliminary results are given in Sections \ref{sec:setup} and \ref{sec:IsoVsConc}, respectively, while in Section \ref{sec:MS}, we introduce the spectral condition we are interested in, relate it to minimal concordance and investigate the matter in the case of surfaces.

\section{Setup and key definitions}\label{sec:setup}

Given $n\geq 2$, we let $X^n$ denote a compact manifold, of dimension equal to $n$, with possibly non-empty boundary. Consistently with our previous work \cite{CarLi21}, we will denote by $\met=\met(X)$ the cone of smooth Riemannian metrics on $X$, and we shall be particularly concerned -- at least initially -- with its topological subspaces defined by binary relations involving its scalar curvature and the mean curvature of its boundary.

For a Riemannian metric $h$ on $X$ we let $\eta=\eta^{(h)}$ denote an outward-pointing unit normal vector field along $\partial X$, take $\II_h$ to be the scalar-valued second fundamental form (with respect to $\eta$) and $H_h$ its trace (that is: the mean curvature of $\partial X$); throughout this article, we adopt the convention that the unit sphere in $\R^3$ has mean-curvature equal to $2$
and we will say -- for a Riemannian metric $h$ on $X$ -- that $(X,h)$ is mean-convex if the mean curvature of such a~manifold is greater or equal than zero, namely if $H_{h}\geq 0$.

In particular, we let $\met_{R>0, H\geq 0}$ (respectively $\met_{R>0, H=0}$) denote the subspace of metrics with positive scalar curvature and mean-convex (respectively, minimal) boundary.

\begin{Definition}\label{eq:Isotopy}
In the setting above, given a topological subspace $\met_{\ast}(X)$ we will say that $h_0, h_1 \in \met_{\ast}(X)$ are isotopic in $\met_{\ast}(X)$ if there exists $\boldsymbol{h}\in C^0(I, \met_{\ast}(X))$ such that $\boldsymbol{h}(0,\cdot)=h_0$ and $\boldsymbol{h}(1,\cdot)=h_1$. (Throughout this article, we set $I:=[0,1]\subset\R$.)
\end{Definition}

 It is straightforward to check that the preceding definition determines an \emph{equivalence relation} in the space $\met_{\ast}(X)$.
One can however consider different equivalence relations within a given subspace of Riemannian metrics on a background manifold $X$. For the purposes of the present article, the following two turn out to be especially significant. In the setting above, let $(\star)$ denote a curvature condition, defined by a binary operator within the set $ \{=,\geq,\leq, >,< \}$, on the mean curvature for $\partial X\times I$ as part of the boundary of the cylinder $X\times I$.

\begin{Definition}\label{def:Concordance}
We will say that $h_0, h_1 \in \met_{\ast}(X)$ are PSC $(\star)$-concordant if there exists a PSC Riemannian metric $g$ on $X\times I$ that
\begin{enumerate}\itemsep=0pt
 \item[(i)] satisfies condition $(\star)$;
 \item[(ii)] restricts\footnote{Here, with slight abuse of terminology, when we refer to the restriction of $g$ to $X\times \{0\}$ we really mean the restriction to the corresponding tangent subspace to $X\times \{0\}$ at each given point, and similarly for $X\times \{1\}$.} to $h_0$ on $X\times \{0\}$, and to $h_1$ along $X\times \{1\}$;
 \item[(iii)] is a product near both $X\times \{0\}$ and $X\times\{1\}$.
\end{enumerate}
In the special case when $(\star)$ is the condition that the mean curvature of $\partial X\times I$ in metric $g$ be identically zero (respectively, be non-negative), we will just say that $h_0, h_1 \in \met_{\ast}(X)$ are PSC min-concordant (respectively, PSC mc-concordant).
\end{Definition}

As we are about to see and discuss, the previous definition (albeit natural) is too rigid for certain purposes, and does not account for the natural interplay between scalar curvature and boundary mean curvature that is apparent since at least \cite{GroLaw80-Spin} (cf.\ \cite{Gro18b, Gro18a}) and played a key role in \cite{CarLi19, CarLi21}. So, here is the amendment we wish to propose:

\begin{Definition}\label{def:WeakConcordance}
We will say that $h_0, h_1 \in \met_{\ast}(X)$ are weakly PSC $(\star)$-concordant if there exists a PSC Riemannian metric $g$ on $X\times I$ that
\begin{enumerate}\itemsep=0pt
 \item[(i)] satisfies condition $(\star)$;
 \item[(ii)] restricts to $h_0$ on $X\times \{0\}$, and to $h_1$ along $X\times \{1\}$;
 \item[(iii)] makes both $X\times \{0\}$ and $X\times\{1\}$ free boundary minimal surfaces.
\end{enumerate}
In the special case when $(\star)$ is the condition that the mean curvature of $\partial X\times I$ in metric $g$ be identically zero (respectively, be non-negative), we will just say that $h_0, h_1 \in \met_{\ast}(X)$ are weakly PSC min-concordant (respectively, weakly PSC mc-concordant).
\end{Definition}

Some comments are appropriate.

\begin{Remark}\label{rem:Compare1}
Item (iii) in Definition \ref{def:Concordance} forces the manifolds $(X,h_0)$, $(X,h_1)$ to have positive scalar curvature and satisfy condition $(\star)$ along $\partial X$; in particular, when $(\star)$ is the requirement $H=0$ (respectively, $H\geq 0$), then $(X,h_0)$ and $(X,h_1)$ will necessarily have positive scalar curvature and minimal (respectively, mean-convex) boundary, so that a posteriori we are anyway defining an equivalence relation in $\met_{R>0,H=0}$ (respectively, $\met_{R>0,H\geq 0}$). This is by no means true in the case of the weaker requirement~(iii) in Definition~\ref{def:WeakConcordance}: to make a basic example, we will see later (as a consequence of Corollary~\ref{cor:IsoConcSpecSurf}) that the flat metric on the unit disk $\mathbb{D}$ is weakly PSC min-concordant to the hemispherical metric on the upper half of~$\mathbb{S}^2$.
\end{Remark}

\begin{Remark}\label{rem:Compare2}
Already in the simpler closed case (so lifting the $(\star)$ requirement on the cylindrical boundary) the resulting notion of weak PSC concordance does place both topological restrictions on $X$ and metric restrictions on $h_0$, $h_1$. On the former front, note that $X\times S^1$ would be forced to support PSC metrics, so, for instance, when $n=2$ we have that $X$ must be a topological sphere for the definition to have any content. The metric restrictions are a bit more subtle, but already apparent from Lemma~\ref{lem:Formulae} below.
\end{Remark}

\begin{Remark}\label{rem:Compare3}
It is not obvious, but nevertheless true that weak PSC $(\star)$-concordance is indeed an equivalence relation in $\met(X)$; specifically, the fact that the relation in question is transitive follows -- as a simple special case -- from our companion work \cite{CarLi23b} (in the spirit of Miao's smoothing theorem, see \cite{Mia02}, which suffices to cover the case when $X$ is a closed manifold). In fact, the same deformation methods allow to show that $h_0, h_1\in\met_{R>0,H=0}$ (respectively, in~$\met_{R>0,H\geq 0}$) are PSC min-concordant (respectively, PSC mc-concordant) if and only if they are weakly PSC min-concordant (resp.\ weakly PSC mc-concordant). In other words, if one restricts a priori to such subspaces, then it is equivalent to require, on top of conditions (i) and (ii), for item (iii) that $X\times\{0\}$ and $X\times\{1\}$ meet $\partial X\times [0,1]$ at a right angle and satisfy any of the (local) geometric conditions in the following hierarchy:
$
\centering
(\text{doubling}) \ \Rightarrow \ (\text{totally geodesic}) \ \Rightarrow \ (\text{minimal}),
$
or even that the metric be a Riemannian product in a neighborhood of the bases.
\end{Remark}

\section[Isotopy in R>0, H=0 vs. minimal concordance]{Isotopy in $\boldsymbol{\met_{R>0, H=0}}$ vs.\ minimal concordance}\label{sec:IsoVsConc}

The following lemma collects some basic slicing formulae concerning warped product metrics on cylinders. (Here, and in the sequel, we will informally employ the word \emph{cylinder} when referring to any smooth manifold of the form $X\times J$ for any interval $J\subset\R$; the cylindrical boundary is by definition $\partial X\times J$, while $X\times\left\{\alpha\right\}$ and $X\times\left\{\beta\right\}$ will be refereed to as the bases of the cylinder, if $J=[\alpha,\beta]$ for $\alpha<\beta\in\R$.)

\begin{Lemma}\label{lem:Formulae}
Let us consider on the manifold $M=X\times J$ a smooth metric of the form
\[
g(x,t)=u(x,t)^2 {\rm d}t^2+h_t(x),
\]
where $u\in C^{\infty}(M)$ and the map $J\ni t \mapsto h_t(x)\in \met(X)$ is also smooth. Then the following equations hold:
\begin{enumerate}\itemsep=0pt
\item[$(1)$] $2$nd fundamental form of the slice $X\times\{t\}$
\[\II_t(x)=(2u(x,t))^{-1}\frac{\rm d}{{\rm d}t}h_t(x);\]
\item[$(2)$] mean curvature of the slice $X\times\{t\}$
\[
H_t(x)=(2u(x,t))^{-1}\tr_{h_t}\frac{\rm d}{{\rm d}t}h_t(x);\]
\item[$(3)$] scalar curvature of the product manifold
\[
R_g=2u(x,t)^{-1}\left(-\lapl_{h_t} u+\frac{1}{2}R_{h_t} u\right)-2u(x,t)^{-1}\frac{\rm d}{{\rm d}t}H_t(x)-(H_t(x))^2-|{\II_t}|^2;
\]
\item[$(4)$] mean curvature of the cylindrical boundary of the product manifold
\[
H_g=(u(x,t))^{-1} (H_h u +\partial_{\eta}u ).
\]
\end{enumerate}
Note that for the first two equations we have considered $X\times\{t\}$ as boundary of $X\times [0,t]$, i.e., we have worked with respect to the unit normal $u^{-1}\partial_t$; $R_{h_t}$ denotes, instead the scalar curvature of the manifold $(X,h_t)$.
\end{Lemma}

\begin{proof}
 All formulae but the last one are part of Lemma A.1 of \cite{LiMan21}. So let us briefly discuss item~(4). Said $\left\{\tau_1, \ldots, \tau_{n-1}\right\}$ a local orthonormal frame for the tangent space to $\partial X$, we can complete it to a local orthonormal frame for $TM$ by adjoining a unit normal $\nu$ to $X\times\{t\}$ and an outward unit normal $\eta$ to the cylindrical boundary. In this proof, we will conveniently write~$\langle \cdot, \cdot\rangle$ for the bilinear form corresponding to the metric $g$; let $D$ denote the associated Levi-Civita connection.

By definition of mean curvature, one has $H_g=\sum_{i=1}^{n-1}\langle D_{\tau_i}\eta, \tau_i\rangle+\langle D_{\nu}\eta,\nu\rangle$; here it is clear that
 (with our notation, as in the statement) the first summand equals $H_h$ so we are just left with studying the second. Around a point $p$, let $\left\{x\right\}$ be local coordinates for $X$, so with indices $i,j\in\left\{1,2,\ldots, n\right\}$, and $\nu=u^{-1}\partial_t$, with $t$ henceforth labelled with the index $0$. All sums over repeated indices are from now understood, as it is customary. Hence
$\langle D_{\nu}\eta,\nu\rangle=u^{-2}\langle D_{\partial_0}\eta, \partial_0\rangle$, and, once conveniently set $\varphi:=\langle D_{\partial _0}\eta, \partial_0\rangle$, we have
\[
\varphi=\partial_0\big(\eta^i\big)\langle \partial_i, \partial_0\rangle + \eta^i\langle D_{\partial_0}\partial_i,\partial_0\rangle =\eta^i\langle D_{\partial_0}\partial_i,\partial_0\rangle,
 \]
 where the second equality relies on the block form of the metric in question. We thus need to compute certain Christoffel symbols for a warped product metric as above: indeed
 \[
\langle D_{\partial_0}\partial_i,\partial_0\rangle=\Gamma_{0i}^0 g_{00}
 \]
 and, in turn,
 \[
\Gamma_{0i}^0=\frac{1}{2}g^{0\ell}(g_{0\ell,i}+g_{i\ell, 0}-g_{0i,\ell})=\frac{1}{2}g^{00}(g_{00,i}+g_{i0, 0}-g_{0i,0})=\frac{1}{2}g^{00}g_{00,i}
 \]
 but $g_{00}=u^2(x,t)$, whence $g_{00,i}=2u\partial_i u$. Combining the previous equations we get
 \[
\varphi=\eta^i \Gamma_{0i}^0 g_{00}= \eta^i \frac{1}{2}g^{00}g_{00,i} g_{00}=\frac{1}{2}\eta^i g_{00,i}=\eta^i u\partial_i u
 \]
 and so, at $(x,t)$, there holds $H_g=H_h + u^{-1}\partial_{\eta}u$, as claimed.
\end{proof}

We can then get back to the relation between isotopy and concordance, in our setting. We remark that the implication given in the following statement, in the special case when the base manifold $X$ is a closed manifold, goes back at least to Gromov--Lawson (see \cite[Lemma 3]{GroLaw80-Class}).

\begin{Proposition}\label{pro:IsotVsConc}
Let $(h_t)_{t\in I}$ be an isotopy connecting $h_0$, $h_1$ within $\met_{R>0, H=0}$ $($respectively, within $\met_{R>0, H\geq 0})$. Then $h_0$, $h_1$ are PSC min-concordant $($respectively, PSC mc-concordant$)$.
\end{Proposition}

\begin{proof}
We claim it suffices to consider, on the product $X\times I$, the class of metrics of the form $g(x,t)=A^2 {\rm d}t^2+h_t(x)$ for a suitable choice of (large) $A>0$. First of all, by reparametrizing the isotopy in question (without renaming) we can ensure that the map $t\mapsto h_t$ is smooth (cf., e.g., \cite[Proposition~2.1]{CabMia18}), $h_{t}=h_0$ for $t\in [0,1/3]$ and $h_{t}=h_1$ for $t\in [2/3,1]$.
That being said, Lemma \ref{lem:Formulae} ensures that $H_g=0$ if the given isotopy is in $\met_{R>0, H= 0}$, or else $H_g\geq 0$ if the given isotopy occurs $\met_{R>0, H\geq 0}$ instead; moreover, for the scalar curvature there holds (say at a~point $(x,t)$) the equation $R_g=R_{h_t}+\Omega(x,t)$ where $|\Omega(x,t)|\leq A^{-2}C$ for a positive constant~$C$ only depending of the $C^2$-norm of $I\ni t\mapsto h_t$, in fact just on $\dot{h}_t$, $\ddot{h}_t$, so it suffices to require $A^2>C\rho^{-1}$ for $\rho:=\inf_{(x,t)}R_{h_t}(x)$ to ensure that the given concordance be PSC (i.e., that the metric $g$ on $X\times I$ has positive scalar curvature). \end{proof}

In particular, we wish to stress some remarkable consequences that come straight from the main result obtained by the authors in \cite{CarLi21}:

\begin{Corollary}\label{cor:Gen3Dconcordance}
Let $X$ be a compact, orientable manifold of dimension $3$. Then any two metrics in $\met_{R>0, H=0}$ $($respectively, within $\met_{R>0, H\geq 0})$ are PSC min-concordant $($respectively, PSC mc-concordant$)$.
\end{Corollary}

What is always known to be true, without any dimensional restriction, is what follows (from~\cite[Section~3]{CarLi19}, \cite[Section 4]{BarHan20}):

\begin{Corollary}\label{cor:Reduction to Doubling}
Let $X$ be a compact, orientable manifold of dimension $n\geq 2$. Then any metric in $\met_{R>0, H=0}$ $(\met_{R>0, H\geq 0})$ is PSC min- $($respectively, mc-$)$ concordant to one with doubling boundary condition.
\end{Corollary}

We shall recall here that a smooth Riemannian metric $g$ on a compact manifold with boundary~$X$ can always be written, near the boundary $\partial X$ in the local form $g=ds^2+h_s(x)$ (with~$s$ denoting the distance, measured with respect to $g$, from the boundary in question), and $g$ is called doubling if
\[
\left[\frac{\partial^{(2\ell+1)} h_s}{\partial s^{(2\ell+1)}}\right]_{s=0}=0 \qquad \text{for all} \quad \ell\in\mathbb{N}
\]
at each point $(x,0)$; note that when $\ell=0$ this is the requirement that the boundary be totally geodesic.

So, in practice, by virtue of the previous corollary for any fixed (compact, orientable) manifold the study of PSC min-concordance $\met_{R>0, H=0}$ may be reduced -- whenever convenient -- to the study of the same relation within the smaller subclass given, e.g., by positive scalar curvature metrics satisfying doubling boundary conditions, hence to the problem of PSC-concordance in presence of an isometry of ``reflection'' type.

The question whether the preceding implication (in Proposition \ref{pro:IsotVsConc}) can actually be reversed, namely whether concordant metric are necessarily isotopic, is well-known (see, e.g., \cite{Ros07,RosSto01} as well as references therein) and is still open in its full generality even in the category of closed (boundaryless, compact) manifolds. Here we wish to pose its natural counterpart in our setting:

\begin{question}\label{q:ConcVsIso}
 Is it always the case that PSC min-concordant metrics are in fact isotopic in the space $\met_{R>0, H=0}$?
\end{question}

Of course one can pose an analogous question under different ``boundary conditions'', such as the useful mean-convexity requirement, or more general constraints involving the mean curvature of the boundaries in play. Exactly as for the closed case, the question above is to be understood ``modulo potential topological obstructions'', so for instance in the (important) special case when the base manifold is simply connected.

\section{Spaces of metrics defined by a spectral stability condition}\label{sec:MS}

In this section we wish to discuss how one can effectively design isotopies (and concordances) of metrics in the case of surfaces, namely when $n=2$. The reader may wish to compare, in terms of results and methods, the discussion below with that presented in \cite{CarWu21} by the first-named author and Wu; the main theorem there provides a general criterion ensuring that certain subspaces of $\met(X)$ be either empty or contractible. We note that this conclusion applies, in particular, for~$\met_{R>0, H\geq 0}$ and $\met_{R>0, H=0}$ (such spaces are not empty if and only if~$X$ is a~disk; of course~$R$ is twice the Gauss curvature while $H$ is the geodesic curvature of the boundary), see \cite[Corollary~1.2]{CarWu21}.

As anticipated in the introduction, we will now introduce a different subspace of $\met(X)$, defined by a spectral condition. Given any Riemannian metric $h$ on $X$, for instance possibly the one obtained by pull-back via an immersion in some ambient manifold, we consider the elliptic eingenvalue problem with oblique boundary conditions:
 \be\label{eq:BVP}
 \begin{cases}
 -\lapl_h u+\frac{1}{2}R_h u =\lambda u & \text{on} \ X, \\
 \partial_\eta u + H_h u =0 & \text{on} \ \partial X ,
 \end{cases}
 \ee
 where $\eta$ denotes the outward-pointing unit vector field along the boundary $\partial X$ and $H_h$ is computed consistently with this choice. Then we set
\begin{align}
 \M:={}& \{h \in \met(X) \colon \lambda_1>0,\ \text{where $\lambda_1$ is the principal eigenvalue for \eqref{eq:BVP}}\}\label{eq:MetSpace}
 \end{align}
and explicitly note that such a principal eigenvalue admits the variational characterization:
\be\label{eq:VarChar}
\lambda_1=\inf_{u\neq 0}\frac{\int_X \big(|\nabla_X u|^2+\frac{1}{2}R_h u^2\big){\rm dvol}_h+\int_{\partial X} H_h u^2 {\rm dvol}^{\partial}_h}{\int_{X}u^2 {\rm dvol}_h},
\ee
where, with slightly unconventional notation, we have denoted by ${\rm dvol}_h$ the Riemannian volume element in metric $h$, and ${\rm dvol}^{\partial}_h$ the $(n-1)$-dimensional Riemannian volume element induced, by the same metric, on the boundary $\partial X$. When we wish to stress the dependence on the background metric, for instance to avoid ambiguities, we will write $\lambda^{(h)}_1$.

Furthermore, we shall remark that condition \eqref{eq:BVP} implies, for $h\in \M$, that one can employ the associated first eigenfunction to design a Riemannian metric on $X\times I$ satisfying peculiar curvature conditions.

\begin{Lemma}\label{lem:WarpFirstEigenf}
In the setting above, let $u=u(x)>0$ denote a first eigenfunction for the oblique eigenvalue problem \eqref{eq:BVP}, and let $\lambda_1$ denote the corresponding eigenvalue. If $h\in \M$, then the metric $g=u(x)^2{\rm d}t^2+h$ on $X\times I$ has positive scalar curvature $($equal to $2\lambda_1)$ and minimal boundary.
\end{Lemma}

\begin{proof}
This is straightforward from the formulae collected in Lemma~\ref{lem:Formulae}, since by construction all slices $X\times \{t\}$ of the cylinder in question are totally geodesic.
\end{proof}

 The introduction of the space $\M$ is justified by the following statement, whose proof is an application of the Schoen--Yau rearrangement trick.

 \begin{Lemma}\label{lem:CharFirstEigenv}
Let $\big(M^{n+1},g\big)$ be a Riemannian manifold with non-negative scalar curvature ${(R_g\geq 0)}$, weakly mean-convex boundary $(H_g\geq 0)$, and such that
\be\label{eq:StrictIneq}
\inf_{M} R_g+\inf_{\partial M} H_g>0.
\ee
If $X$ is a compact, connected, properly embedded, two-sided stable free boundary minimal hypersurface, then the induced metric $h=g_X$ $($on the corresponding tangent bundle to $X)$ belongs to~$\M=\M(X)$. Furthermore, when $n=2$, then $\Sigma$ is a topological disk.
 \end{Lemma}

 \begin{proof}
In the context of the statement, the stability inequality reads
\begin{align*}
\int_{X}|\nabla_{X}u|^2 {\rm dvol}_h\geq{}& \int_{X}\big(|\II_h|^2+\Ric_g(\nu,\nu)\big)u^2 {\rm dvol}_h\\
&+\int_{\partial X}\II_g(\nu,\nu)u^2 {\rm dvol}^{\partial}_h, \qquad \forall u\in C^{\infty}_c(X),
\end{align*}
where $\II_g$ denotes the (scalar-valued) second fundamental form of $\partial M$, $\II_h$ is as above, and $\nu$ is a~choice of the unit normal to $X$ inside $M$. Such an inequality
can be conveniently rearranged~to
\begin{align*}
&\int_{X}\left(|\nabla_{X}u|^2+\frac{1}{2}R_h u\right)u {\rm dvol}_h + \int_{\partial X}H_h u^2\, {\rm dvol}^{\partial}_h \\
 &\qquad\geq \int_{X}\frac{1}{2}\big(\Scal_g+|\II_h|^2\big)u^2 {\rm dvol}_h+\int_{\partial X}(\operatorname{tr}_g\II_g) u^2 {\rm dvol}^{\partial}_h, \qquad \forall u\in C^{\infty}_c(X).
\end{align*}
At this stage, the assumption \eqref{eq:StrictIneq} (noting that $\operatorname{tr}_g\II_g=H_g$) together with the very definitions~\eqref{eq:MetSpace} and \eqref{eq:VarChar} implies the first conclusion. Instead, the second claim (corresponding to~$n=2$) follows from taking $u=1$ in the previous inequality and appealing to the Gauss--Bonnet theorem since the Euler characteristic of $X$ is forced to be strictly positive.
 \end{proof}

We note that the previous statement can actually by upgraded to a characterization:

\begin{Proposition}\label{prop:EquivCharStab}
 Let $X^n$ be a compact manifold with boundary. Then the metrics lying in~$\M(X)$ are all and only those induced by the embedding of $X$ as a stable, free boundary minimal hypersurface in a manifold $\big(M^{n+1},g\big)$ satisfying $R_g\geq 0$, $H_g\geq 0$ and \eqref{eq:StrictIneq}.
\end{Proposition}

\begin{proof}
 One implication has been discussed in Lemma~\ref{lem:CharFirstEigenv}, so let us see the converse. Given $h\in \M(X)$, consider the manifold $M=X\times S^1$ endowed with the metric $g=u^2 {\rm d}t^2+h$ where $u=u(x)>0$ is a first eigenfunction for \eqref{eq:BVP}. As we have seen above, in Lemma~\ref{lem:WarpFirstEigenf}, the metric~$g$ has positive scalar curvature and minimal boundary, so condition \eqref{eq:StrictIneq} is certainly satisfied (thanks to the compactness of $X$). Furthermore, fixed any $t_0\in S^1$ one has that~$X\times \{t_0\}$ is a two-sided totally geodesic hypersurface, meeting the boundary of the ambient metric orthogonally; its stability follows, e.g., from Barta's criterion (cf.\ \cite[Lemma~1.36]{ColMin11}) since the function $\equiv 1$ patently lies in the kernel of its Jacobi operator.
\end{proof}

Next, we prove that $\M(X)$ -- when not empty -- has the simplest possible homotopy type:

\begin{Theorem}\label{thm:Contract}
 Let $X=\mathbb{D}^2$. Then the space $\M(X)$ is contractible.
\end{Theorem}

\begin{proof}
 We will verify that the space $\M(X)$ satisfies the assumptions of the aforementioned Theorem 1.1 in \cite{CarWu21}. That this space is invariant under diffeomorphism is obvious, so our task is rather to check that the condition defining this space of metrics is ``is convex along the fibers'' in the sense of condition (1) therein. So, let $h_1={\rm e}^{2w_1}h_{\ast}$ and $h_2={\rm e}^{2w_2}h_{\ast}$ both belong to $\M(X)$ and we want to prove that $h_\textbf{t}={\rm e}^{2(t_1 w_1+t_2w_2)}h_{\ast}$ also belongs to $\M(X)$ for any non-negative $t_1$,~$t_2$ such that $t_1+t_2=1$; here $h_{\ast}$ simply denotes a background metric on the disk (under no further specifications). In view of \eqref{eq:BVP} and \eqref{eq:VarChar}, we have to find a positive lower bound for the bilinear form
 \[
 B_{\textbf{t}}(u,u)=\int_X \biggl(\bigl|\nabla^{h_{\textbf{t}}}_X u\bigr|_{h_{\textbf{t}}}^2+\frac{1}{2}R_{h_{\textbf{t}}} u^2\biggr){\rm dvol}_{h_{\textbf{t}}} +\int_{\partial X} H_{h_{\textbf{t}}} u^2 {\rm dvol}^{\partial}_{h_{\textbf{t}}}
 \]
 over the Hilbert sphere $\int_X u^2 {\rm dvol}_{h_{\textbf{t}}}=1$. If we exploit the well-known formulae for the conformal change of Gauss curvature and geodesic curvature we can rewrite (with the obvious notation $\ast$ referring to all sorts of geometric quantities computed in metric $h_{\ast}$)
 \begin{align*}
B_{\textbf{t}}(u,u)={}&\int_X \big(|\nabla^{\ast}_X u|_{h_{\ast}}^2+(K_{\ast} -\Delta_{\ast}(t_1 w_1+t_2 w_2))u^2\big) {\rm dvol}_{h_{\ast}}\\
&+\int_{\partial X} (\kappa_{\ast}+\eta_{\ast}(t_1 w_1+t_2 w_2)) u^2 {\rm dvol}^{\partial}_{h_{\ast}},
 \end{align*}
 which in turn one can rearrange as
 \begin{gather*}
B_{\textbf{t}}(u,u)= t_1\bigg(\int_X \big(|\nabla^{\ast}_X u|_{h_{\ast}}^2+(K_{\ast} -\Delta_{\ast}w_1)u^2\big) {\rm dvol}_{h_{\ast}}+\int_{\partial X} \big(\kappa_{\ast}+\eta_{\ast}(w_1) u^2\big){\rm dvol}^{\partial}_{h_{\ast}}\bigg) \\
\hphantom{B_{\textbf{t}}(u,u)=}{} +t_2 \bigg(\int_X \big(|\nabla^{\ast}_X u|_{h_{\ast}}^2+(K_{\ast} -\Delta_{\ast}w_2) u^2\big) {\rm dvol}_{h_{\ast}}+\int_{\partial X} \big(\kappa_{\ast}+\eta_{\ast}(w_2) u^2\big){\rm dvol}^{\partial}_{h_{\ast}}\bigg).
 \end{gather*}
 As a result, if we let $\rho={\rm e}^{-2\overline{w}}$ for $\overline{w}=\max_{i=1,2}\sup_X |w_i|$, we have
 \begin{align*}
B_{\textbf{t}}(u,u)& \geq t_1 \lambda^{(h_1)}_1 \int_{X}u^2 {\rm dvol}_{h_1} + t_2 \lambda^{(h_2)}_1 \int_{X}u^2 {\rm dvol}_{h_2}
 \\
 &= t_1 \lambda^{(h_1)}_1 \int_{X}u^2 {\rm e}^{2w_1} {\rm dvol}_{h_\ast} + t_2 \lambda^{(h_2)}_1 \int_{X}u^2 {\rm e}^{2w_2} {\rm dvol}_{h_\ast} \\
 & \geq \rho \left( t_1 \lambda^{(h_1)}_1 \int_{X}u^2 {\rm dvol}_{h_{\textbf{t}}} + t_2 \lambda^{(h_2)}_1 \int_{X}u^2 {\rm dvol}_{h_{\textbf{t}}}\right)
 \end{align*}
 and so we can bound from below this quantity with the corresponding convex combinations of the bilinear forms for the metrics $h_1$, $h_2$:
 \[
B_{\textbf{t}}(u,u)\geq \rho \big(t_1 \lambda^{(h_1)}_1+t_2 \lambda^{(h_2)}_1\big)
 \geq \rho \min\big\{\lambda^{(h_1)}_1;\lambda^{(h_2)}_1\big\}>0.
 \]
 Thus we obtain the desired conclusion.
\end{proof}

This argument applies, as a special case, when $X$ is a closed surface \big(i.e., $\mathbb{S}^2$, $\mathbb{R}\mathbb{P}^2\big)$ and thus strengthens \cite[Proposition 1.1]{ManSch15}, cf.\ \cite[Theorem 3.4]{RosSto01} and \cite[Section 7]{AkutagawaBotvinnik2002relative}.

One can in fact further refine the isotopies above and employ Moser's trick to ensure ``constancy of the area form''. Here is the relevant statement:

\begin{Lemma}\label{lem:Moser}
Let $X$ be a compact manifold with boundary, and let $\met_\ast(X)$ be a subspace of $\met(X)$ that is invariant under dilations, as well as under diffeomorphisms. Then the following holds: given any isotopy $(h_t)_{t\in I}$ in $\met_\ast(X)$ connecting two metrics having the same volume, meaning that $\int_X {\rm dvol}_{h_0}=\int_X {\rm dvol}_{h_1}$, there exists an isotopy $\big(\tilde{h}_t\big)_{t\in I}$ also lying in $\met_\ast(X)$, connecting the same metrics $($i.e., with the same endpoints$)$ possibly modulo the action of a pull-back through a~diffeomorphism on $h_1$ and such that, in addition, $\frac{\rm d}{{\rm d}t}{\rm dvol}_{\tilde{h}_t}=0$ for all $t\in I$.
\end{Lemma}

\begin{proof}

Fix $\sigma_t>0$ such that $\sigma_0=\sigma_1=1$ and
\[\frac{\rm d}{{\rm d}t} \vol_{\sigma_t h_t}(X)=0, \qquad t\in I; \] it is clear that these numbers are a continuous function of $t\in I$ and are in fact smooth if so is the path $t\mapsto h_t$ (which by regularization we may and shall assume, without loss of generality);
we further note that
\[\int_X \left(\frac12 \tr_{\sigma_t h_t} \frac{\rm d}{{\rm d}t} (\sigma_t h_t)\right) {\rm dvol}_{\sigma_t h_t} = \frac{\rm d}{{\rm d}t} \vol_{\sigma_t h_t}(X)=0.\]
Thus, for any $t\in I$ we can define $f_t\colon X\to\R$ to be the unique (null mean) solution to
\[\begin{cases} \displaystyle \Delta_{\sigma_t h_t} f_t = -\frac12 \tr_{\sigma_t h_t} \frac{\rm d}{{\rm d}t}(\sigma_t h_t) &\text{in}\ X,\vspace{1mm}\\
\displaystyle \frac{\partial f_t}{\partial \eta_t} = 0 &\text{on}\ \partial X,\end{cases}\]
for $\eta_t$ is the outward pointing unit normal vector field along $\partial X$. Consider then the corresponding gradient vector field, $W_t = \nabla_{\sigma_t h_t} f_t$: the Neumann condition imposed on $f_t$ guarantees that $W_t$ is tangential to $\partial X$ along the boundary. Thus, if we let $\Psi=\Psi(x,t)$ to be the integral flow of the vector field $W_t$, then, set $\psi_t(x)=\Psi(x,t)$, $(\psi_t)$ is a smooth family of boundary-preserving diffeomorphisms of $X$, isotopic to the identity, such that
\[\psi_0 = \text{Id}, \qquad \frac{\rm d}{{\rm d}t}{\rm dvol}_{\psi_t^* (\sigma_t h_t) }=0.\]
Indeed, the relevant computation reads as follows:
\begin{align*}
\frac{\rm d}{{\rm d}t}{\rm dvol}_{\psi_t^* (\sigma_t h_t)}={}&\psi_t^* \left[\frac{\rm d}{{\rm d}t} {\rm dvol}_{\sigma_t h_t}+\mathscr{L}_{\dot{\psi}_t}{\rm dvol}_{\sigma_t h_t}\right] \\
={}& \psi_t^* \left[\left(\frac12 \tr_{\sigma_t h_t} \frac{\rm d}{{\rm d}t} (\sigma_t h_t)+\div_{\sigma_t h_t}\dot{\psi_t}\right){\rm dvol}_{\sigma_t h_t}\right]=0,
\end{align*}
where the last equality relies on the fact that, by construction, $\dot{\psi_t}=W_t=\nabla_{\sigma_t h_t}f_t$ and the very definition of $f_t$.

(The desired computation is local, in particular it does not involve integration by parts, so that the boundary plays no role here; thus, one can just follow that given in the proof of \cite[Lemma 1.2]{ManSch15}.)
Hence, if we simply set $\tilde{h}_t=\psi_t^* (\sigma_t h_t)$ the desired conclusion follows.
\end{proof}

\begin{Corollary}\label{cor:MoserMX}
Let $X$ be a compact manifold with boundary, and let $\M(X)$ be defined as above. Given any isotopy $(h_t)_{t\in I}$ in $\M(X)$ connecting two metrics having the same volume, meaning that $\int_X {\rm dvol}_{h_0}=\int_X {\rm dvol}_{h_1}$, there exists an isotopy $\big(\tilde{h}_t\big)_{t\in I}$ also lying in $\M(X)$, connecting the same metrics $($i.e., with the same endpoints$)$ possibly modulo the action of a pull-back $($through a diffeomorphism$)$ on $h_1$ and such that, in addition, $\frac{\rm d}{{\rm d}t}{\rm dvol}_{\tilde{h}_t}=0$ for all $t\in I$.
\end{Corollary}

When dealing with surfaces, one can further specify and strengthen the result.

\begin{Corollary}\label{cor:MoserSurf}
 Let $X=\mathbb{D}^2$. For any $h\in \M(X)$ there exists an isotopy $(h_t)$ satisfying the following properties:
 \begin{enumerate}\itemsep=0pt
 \item[$(1)$] $h_0=h$ and $h_1$ is diffeomorphic to the standard metric on the unit disk in Euclidean $\R^2$;
 \item[$(2)$] $\dot{h}_t=0$ on $[0,1/3]\sqcup [2/3,1]$;
 \item[$(3)$] $\frac{\rm d}{{\rm d}t}{\rm dvol}_{h_t}=0$.
 \end{enumerate}
\end{Corollary}

\begin{proof}
 First of all, by the uniformization theorem we can write $h_0={\rm e}^{2w}h_{\ast}$ where $h_{\ast}$ is homothetic to the standard metric on the unit disk in Euclidean $\R^2$; since patently $h_{\ast}\in \M(X)$ it follows, as a special case, from the argument given in the proof of Theorem \ref{thm:Contract} that the whole segment ${\rm e}^{2tw}h_{\ast}$, $t\in I$ lies within $\M(X)$. Note that $h_{\ast}$ is chosen so to assign $X$ the same volume (in fact: area) as in metric $h=h_0$.

 By reparametrizing we can certainly ensure that condition (2) above is met as well. At this stage we apply to this isotopy, seen as input, the construction presented in Lemma \ref{lem:Moser}.
\end{proof}

The task of properly understanding the space $\M(X)$ when the underlying manifold $X$ has dimension at least three looks rather challenging. One possible approach -- arguably not the only one -- consists in first relating it to other spaces of Riemannian metrics whose topology (specifically: homotopy type) we understand comparatively well. In particular,
in view of the main results in \cite{CarLi19, CarLi21} we wish to ask the following question:

 \begin{question}\label{qu:DoublingSpectralCondition}
 Let $X$ be a compact manifold of dimension equal to three. Is it true that the space $\M(X)$ is either empty or contractible?
 \end{question}

This question is particularly intriguing, especially on the analytic side, because the definition of $\M(X)$ involves a non-local condition; when $n\geq 4$ the answer to the question above should (often) be negative, although the landscape is at the moment still largely unexplored.

We now get back to surfaces, and turn from the design of isotopies to the design of suitable concordances.

\begin{Theorem}\label{thm:IsoConcSpec}
 Let $X$ be a compact manifold with boundary. If two metrics in $\M(X)$ are isotopic, then they are weakly PSC min-concordant.
\end{Theorem}

\begin{proof}
 Let $(h_t)$ denote an isotopy within $\M(X)$, which we can always smoothly reparametrize so that $\dot{h}_t=0$ on $[0,1/3]\sqcup [2/3,1]$. Now we consider on $M=X\times I$ metrics of the form $g=A^2 u^2_t {\rm d}t^2+h_t$ where $u_t>0$ is a first eigenfunction for \eqref{eq:BVP} \big(normalized so that, say, $\int_X u^2_t\, {\rm dvol}_{h_t}=1$\big) and $A>0$ is a constant to be chosen suitably large at a later stage of the proof. (A functional analytic argument, along the lines of \cite[Appendix~A]{ManSch15}, ensures that the map $t\mapsto u_t \in C^{\infty}(X;\R)$ can be chosen to be itself smooth). Note that the free boundary condition holds because of the block structure of the metric~$g$. From the first item of Lemma~\ref{lem:Formulae} the bases of the cylinder are minimal (in fact totally geodesic); from the last item in the same statement and the oblique boundary condition imposed on $u_t$ it follows that the requirement that the cylindrical boundary of $M$ be minimal is also certainly met. Concerning the scalar curvature of~$M$, we look at the third item therein, and thus observe that (cf.\ proof of Proposition~\ref{pro:IsotVsConc}) \smash{$R_g=2\lambda^{(t)}_1+O\big(A^{-2}\big)$} where \smash{$\lambda^{(t)}_1$} is the principal eigenvalue associated to the eigenfunction~$u_t$, $t\in I$. Hence, said $0<\lambda_{\ast}=\inf_{t\in I}\lambda_t$ we take $A$ large enough to guarantee that the remainder term above be less that~$\lambda_{\ast}$.
\end{proof}

In particular, because of Theorem \ref{thm:Contract} (just at the $\pi_0$ level) when $X=\mathbb{D}^2$ any two metrics in~$\M(X)$ are automatically weakly PSC min-concordant. In particular, we wish to spell out what follows:

\begin{Corollary}\label{cor:IsoConcSpecSurf}
 Let $X=\mathbb{D}^2$. Any $h\in \M(X)$ is weakly PSC min-concordant to the standard metric on the unit disk in Euclidean $\R^2$.
\end{Corollary}

It is clear that the previous statement cannot possibly hold, and does not even make sense, in view of Remark \ref{rem:Compare1}, for (strong) PSC min-concordance.

\begin{Remark}[{outward-bending in the sense of Mantoulidis--Schoen}]
For future reference, we discuss here in what terms and how one can modify the construction of (isotopies and) concordances, so to obtain a conclusion similar to those collected in \cite[Lemma~1.3]{ManSch15}, that being the crucial ancillary step before gluing to exterior Schwarzschild in order to estimate the Bartnik mass of minimal triples.

Let us consider, for $X$ a compact manifold with boundary and $\M(X)$ defined as above, an isotopy as produced by Corollary~\ref{cor:MoserMX}, so a path of metrics in $\M(X)$ preserving the volume form at each point and each time. Then we can place on $M:=X\times I$ a Riemannian metric $g$ of the form
\[
g=A^2u^2_t {\rm d}t^2+\big(1+\e t^2\big)h_t
\]
and, for $\e>0$ sufficiently small and $A>0$ sufficiently large (only depending on $(h_t)_{t\in I}$), we claim that the following statements hold:
\begin{enumerate} \itemsep=0pt
 \item[(i)] the scalar curvature of $(M,g)$ is positive;
 \item[(ii)] the mean curvature of the cylindrical boundary of $(M,g)$ is identically equal to zero;
 \item[(iii)] the slice $X\times \{0\}$ is minimal, and all slices $X\times \{t_0 \}$, $t_0\in (0,1]$ are strictly mean-convex.
\end{enumerate}
(Besides, let us also note that all slices $X\times \{t\}$, $t\in [0,1]$ will also satisfy the free boundary condition, namely each of them will meet the ambient boundary orthogonally.)
Given the \emph{Ansatz} above, the computations are exactly as in the closed case (for which we refer the reader to \cite[pp.~5--6]{ManSch15}) except for the study of the cylindrical boundary. In that respect, if we simply set $k_t=\big(1+\e t^2\big)h_t$, $v_t=Au_t$ one can write $g=v^2_t{\rm d}t^2+k_t$ and we are in the setting of Lemma~\ref{lem:Formulae}, hence in particular \smash{$H_g=(v(x,t))^{-1}(H_k v +\partial_{\eta}v)$} for $\eta=\eta^{(k)}$ the outward-pointing normal in metric $g$. Since such a vector is patently orthogonal to $\partial_t$ there holds (with obvious meaning of the symbols) \smash{$\eta^{(k)}=\bigl(1+\e t^2\bigr)^{-1/2}\eta^{(h)}$}, and -- by appealing to the standard formula for the conformal change of mean curvature -- \smash{$H_k=\bigl(1+\e t^2\bigr)^{-1/2}H_h$}. Thus, it follows that \smash{$H_g=\bigl(1+\e t^2\bigr)^{-1/2}(u(x,t))^{-1}\big(H_h u +\partial_{\eta^{(h)}}u\big)=0$}, where the last equality relies on the boundary conditions imposed on the eigenfunction $u$. Therefore, the conclusion follows.
\end{Remark}

\subsection*{Acknowledgements}

The authors wish to express their sincere gratitude to the editors of this special issue for the possibility of contributing with the present article, and to the anonymous referees for their valuable suggestions.
This project has received funding from the European Research Council (ERC) under the European Union’s Horizon 2020 research and innovation programme (grant agreement No. 947923). C.L. was supported by an NSF grant (DMS-2202343).

\pdfbookmark[1]{References}{ref}
\LastPageEnding

\end{document}